\documentclass[reqno,final]{amsart}
\usepackage{chrome}
\hypersetup{pdfpagemode=UseOutlines,colorlinks=true,pdfpagelayout=SinglePage,pdfstartview=FitH,bookmarksopen=true}

\usepackage[color,notref,notcite]{showkeys}

\newcommand{\gee}{TM\oplus T^*M}
\newcommand{\nij}[2]{\mathcal{N}(#1,#2)}

\newcommand{\II}{\boldsymbol{I}}
\newcommand{\JJ}{\boldsymbol{J}} 
\newcommand{\KK}{\boldsymbol{K}} 

\DeclareMathOperator{\torsion}{T^\nabla}

\begin{document}
\selectlanguage{english}

\title{Hyper-complex structures on Courant algebroids}
\author[M.~Sti\'enon]{Mathieu Sti\'enon}
\address{Universit\'e Paris Diderot, Institut de Math\'ematiques de Jussieu (UMR~CNRS~7586), Site Chevaleret, Case~7012, 
75205~Paris~CEDEX~13, France}
\email{stienon@math.jussieu.fr}
\address{Pennsylvania State University, Department of Mathematics, 109 McAllister Building, University Park, PA 16802, United States}
\email{stienon@math.psu.edu}
\begin{abstract}
Hypercomplex structures on Courant algebroids unify holomorphic symplectic structures and usual hypercomplex structures. 
In this note, we prove the equivalence of two characterizations of hypercomplex structures on Courant algebroids, one in terms of  Nijenhuis concomitants and the other in terms of (almost) torsionfree connections for which each of the three complex structures is parallel.
\end{abstract}
\maketitle

A Courant algebroid \cite{LiuWeinsteinXu} consists of a vector bundle $\pi:E\to M$, a nondegenerate symmetric pairing $\ip{}{}$ on the fibers of $\pi$, a bundle map $\rho:E\to TM$ called anchor and an $\RR$-bilinear operation $\db{}{}$ on $\sections{E}$ called Dorfman bracket, which, for all $f\in\cinf{M}$ and $x,y,z\in\sections{E}$ satisfy the relations 
\begin{align}
& \db{x}{(\db{y}{z})}=\db{(\db{x}{y})}{z}+\db{y}{(\db{x}{z})}; \label{elephant} \\ 
& \rho(\db{x}{y})=\lb{\rho(x)}{\rho(y)}; \\ 
& \db{x}{fy}=\big(\rho(x)f\big)y+f(\db{x}{y}); \\ 
& \db{x}{y}+\db{y}{x}=2\DD\ip{x}{y}; \label{mouse} \\ 
& \db{\DD f}{x}=0; \\ 
& \rho(x)\ip{y}{z}=\ip{\db{x}{y}}{z}+\ip{y}{\db{x}{z}}, 
\end{align}
where $\DD:\cinf{M}\to\sections{E}$ is the $\RR$-linear map defined by $\ip{\DD f}{x}=\thalf\rho(x)f$. 

The symmetric part of the Dorfman bracket is given by \eqref{mouse}. The Courant bracket is defined as the skew-symmetric part $\cb{x}{y}=\thalf(\db{x}{y}-\db{y}{x})$ of the Dorfman bracket. Thus we have the relation $\db{x}{y}=\cb{x}{y}+\DD\ip{x}{y}$. 

A standard example is due to T.~Courant~\cite{Courant}. Given a smooth manifold $M$, the vector bundle $\gee\to M$ carries a natural Courant algebroid structure: the anchor is the projection onto the tangent component while the pairing and Dorfman bracket are given, respectively, by 
\[
\ip{X+\xi}{Y+\eta}=\thalf\big(\xi(Y)+\eta(X)\big) \quad\text{and}\quad 
\db{(X+\xi)}{(Y+\eta)}=\lb{X}{Y}+(\ld{X}\eta-\ii{Y}d\xi),
\]
for all $X,Y\in\XX(M)$ and $\xi,\eta\in\OO^1(M)$.

\begin{defn}
An \emph{almost hypercomplex structure} on a Courant algebroid $(E,\rho,\ip{}{},\db{}{})$ is a triple $(\II,\JJ,\KK)$ of endomorphisms of the vector bundle $E$, i.e. vector bundle maps over $\id_M:M\to M$, which are orthogonal transformations w.r.t. the pairing $\ip{}{}$ and satisfy the quaternionic relations 
\beq{deer} \II^2=\JJ^2=\KK^2=\II\JJ\KK=-1 .\eeq
\end{defn}

Obviously, if $(\II,\JJ,\KK)$ is an almost hypercomplex structure, then so are $(\KK,\II,\JJ)$ and $(\JJ,\KK,\II)$. 

Let $(E\to M,\rho,\ip{}{},\db{}{})$ be a Courant algebroid.
Given two endomorphisms $F$ and $G$ of the vector bundle $E$, the relation 
\begin{equation} \begin{split} 
\nij{F}{G}(X,Y)&=\db{FX}{GY}-F(\db{X}{GY})-G(\db{FX}{Y})+FG(\db{X}{Y}) \\ 
 &\quad +\db{GX}{FY}-G(\db{X}{FY})-F(\db{GX}{Y})+GF(\db{X}{Y})
,\end{split} \label{shark} \end{equation}
where $X,Y\in\sections{E}$, defines a (2,1)-tensor 
$\nij{F}{G}:E\otimes E\to E$ 
called Nijenhuis concomitant. 
Obviously, $\nij{F}{G}=\nij{G}{F}$. 

\begin{lem}
If $(\II,\JJ,\KK)$ is an almost hypercomplex structure on a Courant algebroid $E$, then 
$\nij{\II}{\JJ}(X,Y)+\nij{\II}{\JJ}(Y,X)=0$ for all $X,Y\in\sections{E}$.
\end{lem}

\begin{defn}
A \emph{hypercomplex structure} on a Courant algebroid $E$ is an almost hypercomplex structure $(\II,\JJ,\KK)$ such that 
the six Nijenhuis concomitants $\nij{\II}{\II}$, $\nij{\JJ}{\JJ}$, $\nij{\KK}{\KK}$, $\nij{\II}{\JJ}$, $\nij{\JJ}{\KK}$ and  $\nij{\KK}{\II}$ vanish.
\end{defn}

\begin{rmk} 
Let $(E\to M,\rho,\ip{}{},\db{}{})$ be a Courant algebroid and let $\II$ and $\JJ$ be two endomorphisms of $E$ such that:  
$\II^2=\JJ^2=-1$; $\II$ and $\JJ$ anticommute; both $\II$ and $\JJ$ are orthogonal w.r.t. the pairing $\ip{}{}$; and the three Nijenhuis concomitants $\nij{\II}{\II}$, $\nij{\JJ}{\JJ}$ and $\nij{\II}{\JJ}$ vanish. 
Then it is easy to check that the triple $(\II,\JJ,\II\JJ)$ is a hypercomplex structure on the Courant algebroid.
This is the way Bredthauer originally defined hypercomplex structures in \cite{Bredthauer}. See also \cite{Ezhuthachan}.
\end{rmk}

For any $f\in\cinf{M}$ and $X,Y\in\sections{E}$, let 
\beq{tortue} \Delta_f(X,Y)=\ip{X}{Y}\DD f+\ip{\II X}{Y}\II\DD f+\ip{\JJ X}{Y}\JJ\DD f+\ip{\KK X}{Y}\KK\DD f .\eeq
It is clear that 
\[
\Delta_f(X,\II Y)=\II\Delta_f(X,Y) ,\quad 
\Delta_f(X,\JJ Y)=\JJ\Delta_f(X,Y) ,\quad
\Delta_f(X,\KK Y)=\KK\Delta_f(X,Y) 
\]
and 
\[ \Delta_f(X,Y)+\Delta_f(Y,X)=2\ip{X}{Y}\DD f .\]

\begin{defn}
Let $(\II,\JJ,\KK)$ be an almost hypercomplex structure on a Courant algebroid $(E\to M,\rho,\ip{}{},\db{}{})$. 
A \emph{hypercomplex connection} is an $\RR$-bilinear map 
\[ \sections{E}\otimes\sections{E}\to\sections{E}:(X,Y)\mapsto\nabla_X Y \] 
such that, for all $f\in\cinf{M}$ and $X,Y\in\sections{E}$, we have 
\beq{escargot} \nabla_{fX} Y=f\nabla_X Y \eeq
and 
\beq{limace} \nabla_X(fY)=\big(\rho(X)f\big)Y+f(\nabla_X Y)-\Delta_f(X,Y) .\eeq
Its torsion $\torsion:\sections{E}\wedge\sections{E}\to\sections{E}$ is given by 
\beq{chenille} \torsion(X,Y)=\nabla_X Y-\nabla_Y X-\cb{X}{Y} .\eeq 
\end{defn}

\begin{rmk}
If $L$ is an isotropic subbundle of $E$ stable under $\II$, $\JJ$ and $\KK$, then a hypercomplex connection on $E$ induces a usual $L$-connection on $L$. 
\end{rmk}

The purpose of this note is to establish the following result.

\begin{thm}\label{lion} 
Let $(\II,\JJ,\KK)$ be an almost hypercomplex structure on a Courant algebroid $E$. 
The following assertions are equivalent. 
\begin{enumerate}
\item $\nij{\II}{\JJ}=\nij{\JJ}{\JJ}=0$
\item $\nij{\II}{\JJ}=0$
\item The triple $(\II,\JJ,\KK)$ is a hypercomplex structure, i.e. all six Nijenhuis concomitants 
vanish.
\item There exists a hypercomplex connection $\nabla$ satisfying  
\beq{donkey} \nabla\II=\nabla\JJ=\nabla\KK=0 \eeq and, for all $X,Y\in\sections{E}$, \beq{monkey} \torsion(X,Y)=\II\DD\ip{X}{\II Y}+\JJ\DD\ip{X}{\JJ Y}+\KK\DD\ip{X}{\KK Y} .\eeq  
\item There exists a hypercomplex connection satisfying \eqref{donkey} and \eqref{monkey}; it is unique and given by 
\[ \nabla_X Y=-\thalf \KK\big(\db{\JJ Y}{\II X}-\JJ(\db{Y}{\II X}) 
-\II(\db{\JJ Y}{X})+\JJ\II(\db{Y}{X})\big) .\]
\end{enumerate}
\end{thm}

The remainder of this note is devoted to the proof of this theorem.
Straightforward computations lead to the first two lemmas below, of which the former is a generalization of Theorem~1.1 in \cite{YanoAko1}.

\begin{lem}
Given an almost hypercomplex structure $(\II,\JJ,\KK)$, the relation 
\beq{bird} \nabla_X Y=-\thalf \KK\big(\db{\JJ Y}{\II X}-\JJ(\db{Y}{\II X}) 
-\II(\db{\JJ Y}{X})+\JJ\II(\db{Y}{X})\big) \eeq
defines a hypercomplex connection.
Permuting $\II$, $\JJ$ and $\KK$ cyclically in \eqref{bird}, we obtain two other hypercomplex connections: 
\begin{gather}
\label{oiseau} \nabla'_X Y=-\thalf \II\big(\db{\KK Y}{\JJ X}-\KK(\db{Y}{\JJ X}) 
-\JJ(\db{\KK Y}{X})+\KK\JJ(\db{Y}{X})\big) \\ 
\label{vogel} \nabla''_X Y=-\thalf \JJ\big(\db{\II Y}{\KK X}-\II(\db{Y}{\KK X}) 
-\KK(\db{\II Y}{X})+\II\KK(\db{Y}{X})\big) . 
\end{gather}
\end{lem}

\begin{lem}
Given an almost hypercomplex structure $(\II,\JJ,\KK)$, the hypercomplex connection \eqref{bird} satisfies
\begin{gather}
\nabla_X \JJ=0, \label{dog} \\ 
(\nabla_X\II)Y=\thalf\KK\nij{\II}{\JJ}(X,\II Y)+\thalf\JJ\nij{\II}{\JJ}(X,Y), \label{cat}
\end{gather} 
and 
\begin{multline}\label{horse}
\db{X}{Y}+\thalf\KK\nij{\II}{\JJ}(X,Y)=\nabla_X Y-\nabla_Y X+\DD\ip{X}{Y} \\ 
-\big(\II\DD\ip{X}{\II Y}+\JJ\DD\ip{X}{\JJ Y}+\KK\DD\ip{X}{\KK Y}\big) .
\end{multline}
\end{lem}

\begin{cor}
\label{cow}
Let $(\II,\JJ,\KK)$ be an almost hypercomplex structure on a Courant algebroid $E$. 
If $\nij{\II}{\JJ}=0$, then the hypercomplex connection \eqref{bird} satisfies \eqref{donkey} and~\eqref{monkey}.
\end{cor}

\begin{proof}
We always have $\nabla\JJ=0$ by \eqref{dog}. Since $\nij{\II}{\JJ}=0$, \eqref{cat} implies that $\nabla\II=0$. 
And it follows from $\KK=\II\JJ$ that 
\[ \nabla_X\KK=\db{(\nabla_X\II)}{\JJ}+\db{\II}{(\nabla_X\JJ)}=0 .\]
Thus \eqref{donkey} is proved and \eqref{monkey} follows immediately from \eqref{horse} and the relation 
$\db{x}{y}=\cb{x}{y}+\DD\ip{x}{y}$.
\end{proof}

\begin{lem}
\label{sheep}
Given an almost hypercomplex structure $(\II,\JJ,\KK)$, there exists at most one hypercomplex connection satisfying \eqref{donkey} and~\eqref{monkey}.
\end{lem}

\begin{proof}
Assume there exist two such hypercomplex connections $\nabla^1,\nabla^2$. 
Let \[ \Xi(X,Y)=\nabla^2_X Y-\nabla^1_X Y .\] 
It follows from \eqref{donkey} that 
\[ \Xi(X,\II Y)=\II \Xi(X,Y), \quad \Xi(X,\JJ Y)=\JJ \Xi(X,Y), \quad \Xi(X,\KK Y)=\KK \Xi(X,Y) \] 
and from \eqref{monkey} that $\Xi(X,Y)=\Xi(Y,X)$. 
Therefore
\begin{multline*} 
\KK\Xi(X,X) = \II\JJ\Xi(X,X) = \II\Xi(X,\JJ X) = \II\Xi(\JJ X,X) = \Xi(\JJ X,\II X) \\ 
= \Xi(\II X,\JJ X) = \JJ\Xi(\II X,X) = \JJ\Xi(X,\II X) = \JJ\II\Xi(X,X) = -\KK\Xi(X,X) 
.\end{multline*} 
Hence $\Xi(X,X)=0$ for all $X\in\sections{E}$ and, consequently, 
\[ \Xi(X,Y)=\thalf\big(\Xi(X+Y,X+Y)-\Xi(X,X)-\Xi(Y,Y)\big)=0 \] 
for all $X,Y\in\sections{E}$. 
\end{proof}

\begin{lem} 
\label{squirrel}
Given an almost hypercomplex structure $(\II,\JJ,\KK)$, if there exists a hypercomplex connection 
satisfying \eqref{donkey} and \eqref{monkey}, then $\nij{\II}{\JJ}=0$.
\end{lem}

\begin{proof}
From \eqref{monkey}, it follows that 
\begin{multline*} \db{X}{Y}=\nabla_X Y-\nabla_Y X+\DD\ip{X}{Y} 
-\big(\II\DD\ip{X}{\II Y}+\JJ\DD\ip{X}{\JJ Y}+\KK\DD\ip{X}{\KK Y}\big) .\end{multline*}
This relation can be used to evaluate each of the terms of $\nij{\II}{\JJ}$. 
It follows from \eqref{donkey}, the quaternionic relations \eqref{deer}, and the orthogonality of the endomorphisms $\II$, $\JJ$ and $\KK$ w.r.t. the pairing that $\nij{\II}{\JJ}$ vanishes.
\end{proof}

Together, Lemma~\ref{squirrel}, Corollary~\ref{sheep} and Lemma~\ref{cow} imply the following 

\begin{prop}
\label{snake}
Given an almost hypercomplex structure $(\II,\JJ,\KK)$ on a Courant algebroid $E$, 
there exists a hypercomplex connection 
satisfying \eqref{donkey} and \eqref{monkey} 
if and only if $\nij{\II}{\JJ}=0$. And in that case, it coincides with all three hypercomplex connections 
given by \eqref{bird}, \eqref{oiseau} and \eqref{vogel}.
\end{prop}

\begin{prop}
\label{hyena}
Let $(\II,\JJ,\KK)$ be an almost hypercomplex structure on a Courant algebroid. The following assertions are equivalent: 
\begin{enumerate}
\item $\nij{\II}{\II}=\nij{\JJ}{\JJ}=0$;
\item $\nij{\II}{\JJ}=0$; 
\item $\nij{\II}{\II}=\nij{\JJ}{\JJ}=\nij{\KK}{\KK}=\nij{\II}{\JJ}=\nij{\JJ}{\KK}=\nij{\KK}{\II}=0$. 
\end{enumerate}
\end{prop}

\begin{proof}
(i)$\Rightarrow$(ii) The proof is a lengthy computation similar to that of \cite[Theorem~3.1]{YanoAko2}. It is omitted.  
\quad (ii)$\Rightarrow$(iii) 
For any pair of elements $P,Q$ in $\{\II,\JJ,\KK\}$, 
we can evaluate the Nijenhuis concomitant 
\beq{dolphin} \begin{split} 
\nij{P}{Q}(X,Y)&=\db{PX}{QY}-P(\db{X}{QY})-Q(\db{PX}{Y})+PQ(\db{X}{Y}) \\ 
&\quad +\db{QX}{PY}-Q(\db{X}{PY})-P(\db{QX}{Y})+QP(\db{X}{Y})
\end{split} \eeq
by successively making use of: 
\emph{primo} relation \eqref{horse} to get rid of all the Dorfman brackets in the r.h.s. of \eqref{dolphin}; 
\emph{secondo} \eqref{donkey} and the quaternionic relations \eqref{deer} to cancel all terms involving $\nabla$; 
\emph{terzo} \eqref{deer} and the orthogonality of $\II$, $\JJ$ and $\KK$ w.r.t. the pairing to cancel all remaining terms. 
\quad (iii)$\Rightarrow$(i) 
This is trivial.
\end{proof}

Theorem~\ref{lion} immediately follows from Propositions~\ref{snake} and~\ref{hyena}.

\begin{ex}
Let $i$, $j$, $k$ be three almost complex structures on a smooth manifold $X$. 
The triple 
\[ \II = \begin{pmatrix} -i & 0 \\ 0 & i^* \end{pmatrix} ,\quad
\JJ = \begin{pmatrix} -j & 0 \\ 0 & j^* \end{pmatrix} ,\quad
\KK = \begin{pmatrix} -k & 0 \\ 0 & k^* \end{pmatrix} \] 
is a hypercomplex structure on $TX\oplus T^*X$ 
if and only if the triple $i$, $j$, $k$ is hypercomplex in the classical sense (see \cite{Obata}).
\end{ex}

\begin{ex}
Let $j$ be an almost complex structure on a smooth manifold $X$ and let $\omega_1$ and $\omega_2$ be two nondegenerate 2-forms on $X$. 
The triple 
\[ \II = \begin{pmatrix} 0 & \omega_2\inv \\ -\omega_2 & 0 \end{pmatrix} ,\quad
\JJ = \begin{pmatrix} -j & 0 \\ 0 & j^* \end{pmatrix} ,\quad
\KK = \begin{pmatrix} 0 & \omega_1\inv \\ -\omega_1 & 0 \end{pmatrix} \] 
is hypercomplex on $TX\oplus T^*X$ if and only if 
$\omega_1+i\omega_2\in\OO^2_{\CC}(X)$ is a holomorphic symplectic structure on $X$.
Theorem~\ref{lion} has interesting consequences in this case, which we will discuss somewhere else.
\end{ex}

\bibliographystyle{amsplain}
\bibliography{citron}
\end{document}